\documentclass{article}

\usepackage{amsfonts, amsmath, amssymb, amsthm, amscd}
\usepackage{ascmac}

\usepackage{verbatim}

\usepackage{graphics}
\usepackage[all,2cell]{xypic}
\usepackage{mathptmx}

\usepackage{helvet}
\usepackage{courier}        
\usepackage{type1cm} 

\objectmargin+{1mm}
\labelmargin+{0.8mm}
\SelectTips{cm}{12}

\theoremstyle{plain}  
\newtheorem{thm}{Theorem}
\newtheorem{con}[thm]{Conjecture}
\newtheorem{cor}[thm]{Corollary}
\newtheorem{lem}[thm]{Lemma}

\theoremstyle{definition}

\newtheorem{df}[thm]{Definition}
\newtheorem{ex}[thm]{Example}

\newtheorem{conv}[thm]{Conventions}

\theoremstyle{remark}

\DeclareMathOperator{\cB}{\mathcal{B}}
\DeclareMathOperator{\cC}{\mathcal{C}}
\DeclareMathOperator{\calD}{\mathcal{D}}
\DeclareMathOperator{\cE}{\mathcal{E}}

\DeclareMathOperator{\cI}{\mathcal{I}}

\DeclareMathOperator{\calK}{\mathcal{K}}

\DeclareMathOperator{\cM}{\mathcal{M}}

\DeclareMathOperator{\cP}{\mathcal{P}}

\DeclareMathOperator{\cY}{\mathcal{Y}}

\DeclareMathOperator{\bQ}{\mathbf{Q}}

\DeclareMathOperator{\bbH}{\mathbb{H}}

\DeclareMathOperator{\bbK}{\mathbb{K}}

\DeclareMathOperator{\ff}{\mathfrak{f}}

\DeclareMathOperator{\fJ}{\mathfrak{J}}

\DeclareMathOperator{\fp}{\mathfrak{p}}

\DeclareMathOperator{\fs}{\mathfrak{s}}

\DeclareMathOperator{\fY}{\mathfrak{Y}}

\UseAllTwocells

\def\sn{\smallskip\noindent}
\def\mn{\medskip\noindent}

\def\enumidef{\renewcommand{\labelenumi}{$\mathrm{(\Roman{enumi})}$}}
\def\enumiidef{\renewcommand{\labelenumii}{$\mathrm{(\alph{enumii})}$}}

\newcommand{\ab}{\operatorname{ab}}

\newcommand{\CH}{\operatorname{CH}}
\newcommand{\Ch}{\operatorname{\bf Ch}}
\newcommand{\codim}{\operatorname{codim}}

\newcommand{\Cone}{\operatorname{Cone}}
\newcommand{\colim}{\operatorname{colim}}
\newcommand{\coker}{\operatorname{Coker}}

\newcommand{\depth}{\operatorname{depth}}

\newcommand{\Homo}{\operatorname{H}}

\newcommand{\height}{\operatorname{ht}}

\newcommand{\id}{\operatorname{id}}

\newcommand{\isoto}{\overset{\scriptstyle{\sim}}{\to}}

\newcommand{\Ker}{\operatorname{Ker}}

\newcommand{\Mor}{\operatorname{Mor}}

\newcommand{\Nat}{\operatorname{Nat}}

\newcommand{\Ob}{\operatorname{Ob}}
\newcommand{\onto}[1]{\stackrel{#1}{\to}}

\newcommand{\pd}{\operatorname{Projdim}}

\newcommand{\qis}{\operatorname{qis}}

\newcommand{\rdef}{\twoheadrightarrow}
\newcommand{\red}{\operatorname{red}}
\newcommand{\rinc}{\hookrightarrow}
\newcommand{\rinf}{\rightarrowtail}

\newcommand{\simp}{\operatorname{simp}}
\newcommand{\Spec}{\operatorname{Spec}}

\newcommand{\ssm}{\smallsetminus}

\newcommand{\Typ}{\operatorname{Typ}}

\newcommand{\UD}{\operatorname{UD}}
\newcommand{\UT}{\operatorname{UT}}

\title{Local Gersten's conjecture for regular system of parameters}
\date{}

\author{Satoshi Mochizuki}

\begin{document}

\maketitle

\section*{Introduction}

The purpose of this article is to give a proof of 
Gersten's conjecture~\ref{con:Gersten's conjecture} below for certain cases. 
To recall the precise statement of Gersten's conjecture, 
we start by recalling the following result. 
Let $X$ be a smooth variety over a field. 
Then we have the canonical isomorphisms
\begin{equation}
\label{eq:Bloch formula}
\CH^n(X)\simeq\Homo_{\operatorname{Zar}}^n(X,\calK_n)
\end{equation}
where $\calK_n$ is the Zariski sheafication of 
$K$-presheaf $U\mapsto K_n(U)$ and $\CH^n(X)$ is the $n$-th Chow group of $X$. 
For $n=0$, $1$, the result are classical and for $n=2$, 
it was given by Spencer Bloch in \cite{Blo74} 
by utilizing the second universal Chern class map in \cite{Gro68}. 
For general $n$, the proof was given by Quillen in \cite{Qui73} as a consequence of Gersten's conjecture. 
As giving an explanation below, 
in \cite{Ger73}, 
Gersten made an interpretation of the isomorphism 
$\mathrm{(\ref{eq:Bloch formula})}$ 
as the isomorphism of $E_2$-terms of the following two spectral sequences. 
For a noetherian scheme $X$, we write $\cM_X$ for the category of 
coherent sheaves on $X$. 
There is a filtration
\begin{equation}
\label{eq:codimension filtration}
0\subset\cdots\subset\cM_X^2\subset\cM_X^1\subset\cM_X^0=\cM_X
\end{equation}
by the Serre subcategories $\cM_X^p$ of those coherent sheaves whose 
support has codimension$\geq p$. 
Then Quillen showed in \cite{Qui73} that 
for a noetherian scheme $X$, the filtration $\mathrm{(\ref{eq:codimension filtration})}$ 
induces the strongly convergent spectral sequence
\begin{equation}
\label{eq:Quillen spectral sequence}
E_1^{p,q}(X)=\bigoplus_{x\in X^p}K_{-p-q}(k(x))\Rightarrow G_{-p-q}(X),
\end{equation}
where $X^p$ is the set of points codimension $p$ in $X$. 
Moreover if $X$ is regular separated, then we have the canonical isomorphism
\begin{equation} 
\label{eq:E_2 and Chow group}
E_2^{p,-p}(X)\simeq \CH^p(X).
\end{equation}
On the other hand, Brown-Gersten \cite{BG73} proved and 
Thomason \cite{TT90} improved that 
for a noetherian scheme of finite Krull dimension $X$, we have the canonical 
equivalence 
\begin{equation}
\label{eq:Zariski descent}
\bbK(X)\simeq\bbH_{\operatorname{Zar}}^{\bullet}(X;\bbK(-)).
\end{equation}
In particular, there exists the strongly convergent spectral sequence
\begin{equation}
\label{eq:Brown-Gersten spectral seqnece}
{}'E_2^{p,q}(X)=\Homo_{\operatorname{Zar}}^p(X,\tilde{\bbK}_q)\Rightarrow \bbK_{q-p}(X)
\end{equation}
where $\tilde{\bbK}_q$ is the Zariski sheafication of the presheaf 
$U\mapsto \pi_q\bbK(U)$ of non-connective $\bbK$-theory. 
Gersten demonstrated 
the following three conditions are equivalent in \cite{Ger73}. 
For simplicity, for a commutative noetherian ring $A$ with $1$, 
we write $\cM_A$ and $\cM_A^p$ for $\cM_{\Spec A}$ and $\cM^p_{\Spec A}$ 
respectively:

\begin{enumerate}
\item
For any regular separated noetherian scheme $X$, 
we have the canonical isomorphism between $E_2$-terms of Quillen and 
Brown-Gersten-Thomason spectral sequences
$$E_2^{p,q}(X)\simeq{}' E^{p,-q}_2(X).$$
In particular the isomorphism $\mathrm{(\ref{eq:Bloch formula})}$ holds for $X$ and 
any non-negative integer $n$. 

\item
For any commutative regular local ring $R$ of dimension $n$, 
$E_1$-term of Quillen spectral sequence 
$\mathrm{(\ref{eq:Quillen spectral sequence})}$ 
for $\Spec R$ yields an exact sequence
\begin{equation}
\label{eq:Gersten complex}
0\to K_n(R)\to K_n(\operatorname{frac}(R))\to\bigoplus_{\height \fp=1}K_{n-1}(k(\fp)) \to \bigoplus_{\height \fp=2}K_{n-2}(k(\fp))\to \cdots.
\end{equation}

\item
For any commutative regular local ring $R$ and natural number $1\leq p\leq \dim R$, the canonical inclusion $\cM_R^p\rinc \cM_R^{p-1}$ induces the zero map 
on $K$-theory 
$K(\cM_R^p)\to K(\cM_R^{p-1})$ 
where $K(\cM_R^p)$ denotes the $K$-theory of the abelian category $\cM_R^p$.
\end{enumerate}

\sn
Here is Gersten's conjecture:

\begin{con}[\bf Gersten's conjecture]
\label{con:Gersten's conjecture}
The conditions above are true for any commutative regular local ring.
\end{con}

\sn
Corollaries of our main theorem~\ref{main theorem} are the following:

\begin{cor}
\label{main theorem cor 1}
For any commutative regular local ring $R$ of Krull dimension $d$, 
the inclusion functor $\cM^d_R\rinc \cM^{d-1}_R$ induces the zero map 
on $K$-theory 
$K(\cM_R^d)\to K(\cM_R^{d-1})$ 
where $K(\cM_R^p)$ denotes the $K$-theory of the abelian category $\cM_R^p$ 
for $p=d-1$, $d$. 
In particular if $d=1$, namely if $R$ is a commutative discrete valuation ring, we obtain 
Gersten's conjecture for $R$.
\end{cor}

\sn
A proof of Corollary~\ref{main theorem cor 1} 
will be given after Corollary~\ref{cor:LGC for rsp} of 
our main theorem~\ref{main theorem}. 
By virtue of the result by Gillet and Levine \cite[Corollary 6]{GL87} 
(see below), 
Corollary~\ref{main theorem cor 1} implies the following.

\begin{cor}
\label{main theorem cor 2}
Let $R$ be a commutative regular local and smooth over a commutative discrete 
valuation ring $S$. 
Then Gersten's conjecture for $R$ is true. 
\qed
\end{cor}

\subsection*{Historical note}
\label{subsec:Historical note}

The conjecture has been proved in the following cases.

\subsubsection*{The case of general dimension}

\begin{enumerate}
\item
If $A$ is of equi-characteristic, then Gersten's conjecture for $A$ is true. 
We refer to \cite{Qui73} for special cases, and the general cases \cite{Pan03} 
can be deduced from limit argument and 
Popescu's general N\'eron desingularization \cite{Pop86} 
(for the case of commutative discrete valuation rings, it was first proved by 
Sherman \cite{She78}).

\item
If $A$ is smooth over some commutative discrete valuation ring $S$ 
and satisfies some condition, then Gersten's conjecture for $A$ is true 
\cite{Blo86}.

\item
If $A$ is smooth over some commutative discrete valuation ring $S$ 
and if we accept Gersten's conjecture for $S$, then 
Gersten's conjecture for $A$ is true \cite{GL87} (see also \cite{RS90}).
\end{enumerate}

\subsubsection*{The case of that $A$ is a commutative discrete valuation ring}
\label{subsubsec:The case of that A is a commutative discrete valuation ring}

\begin{enumerate}
\item
For the case of $n=0$, $1$ are classical and $n=2$ was proved by 
Dennis and Stein, announced in \cite{DS72} and proved in \cite{DS75}.

\item
If the residue field is algebraic over a finite field, 
then Gersten's conjecture for $A$ is true. 
We refer to \cite{Ger73} for the case of that its residue field is 
a finite field, the proof of general cases \cite{She82} is 
improved from that of special cases by using Swan's result \cite{Swa63}, 
the universal property of algebraic $K$-theory \cite{Hil81} and limit 
argument.
\end{enumerate}

\subsubsection*{The case for Grothendieck groups}

Gersten's conjecture for Grothendieck groups has 
several equivalent forms. 
Namely we can show the following three conditions are equivalent. 
See \cite{CF68}, \cite{Lev85}, \cite{Dut93} and \cite{Dut95}.

\begin{enumerate}
\item[]
{\bf (Gersten's conjecture for Grothendieck groups).}\ 
For any commutative regular local ring $R$ and any natural number 
$1\leq p\leq \dim R$, the canonical inclusion 
$\cM_R^p\rinc \cM_R^{p-1}$ induces the zero map on Grothendieck groups 
$K_0(\cM_R^p)\to K_0(\cM_R^{p-1})$.
\item[]
{\bf (Generator conjecture).}\ 
For any commutative regular local ring $R$ and any natural number 
$0\leq p\leq \dim R$, 
the Grothendieck group $K_0(\cM_R^p)$ is generated by 
cyclic modules $R/(f_1,\cdots,f_p)$ where the 
sequence $f_1,\cdots,f_p$ forms an $R$-regular sequence.
\item[]
{\bf (Claborn and Fossum conjecture).}\ 
For any commutative regular local ring $R$, 
the Chow homology group $\CH_k(\Spec R)$ is trivial for any $k<\dim R$. 
\end{enumerate}

\subsubsection*{The case for singular varieties}
\label{subsubsec:The case for singular varieties}

Gersten's conjecture for non-regular rings is in general false in the 
literal sense of the word and several appropriate modified version are 
studied by many authors. 
See \cite{DHM85}, \cite{Smo87}, \cite{Lev88}, \cite{Bal09}, 
\cite{HM10}, \cite{Moc13a}, \cite{Mor15} and \cite{KM16}.

\subsubsection*{Other cohomology theories}

\begin{enumerate}
\item
For torsion coefficient $K$-theory, Gersten's conjecture 
for a commutative discrete valuation ring is true \cite{Gil86}, \cite{GL00}.

\item
For Gersten's conjecture for Milnor $K$-theory, 
see \cite{Ker09} and \cite{Dah15}.

\item
For Gersten's conjecture for Witt groups, 
see \cite{Par82}, \cite{OP99}, \cite{BW02} and \cite{BGPW03}.

\item
For an analogue of Gersten's conjecture for bivariant $K$-theory, 
see \cite{Wal00}.

\item
For an analogue of Gersten's conjecture by Quillen in 
\cite{Qui73}, he introduced a strengthening of 
Noether normalization theorem. 
There are several variants of Noether normalization theorem in 
\cite{Oja80}, \cite{GS88}, \cite{Gab94} and \cite{Wal98} and 
by utilizing them, there exists 
Gersten's conjecture type theorem for universal exactness \cite{Gra85} 
for Cousin complexes \cite{Har66} of certain cohomology theories. 
For axiomatic approaches of these topics, see \cite{BO74} and 
\cite{C-THK97}. 
See also \cite{Gab93}, \cite{C-T95}, \cite{Lev08} and \cite{Lev08}.

\item
For an analogue of Gersten's conjecture for 
Hochschild coniveau spectral sequences, see \cite{BW16}.

\item
For an analogue of Gersten's conjecture for 
infinitesimal theory, see \cite{DHY15}.

\item
For Gersten's complexes for homotopy invariant Zariski sheaves with transfers, 
see \cite{Voe00}, \cite[Lecture 24]{MVW06} and \cite{SZ03}. 
For injectivity result for pseudo pretheory, see \cite{FS02}.
\end{enumerate}

\subsubsection*{Logical connection with other conjectures}
\label{subsubsec:Logical connection with other conjectures}

\begin{enumerate}
\item
Some conjectures imply that Gersten's conjecture for a commutative 
discrete valuation ring is true \cite{She89}.

\item
Parshin conjecture in \cite{Bei84} implies Gersten's conjecture 
of motivic cohomology for a localization of smooth varieties over 
a Dedekind ring, see \cite{Gei04}. 
It is known that Tate-Beilinson conjecture 
(see \cite{Tat65}, \cite{Bei87} and \cite{Tat94}) implies 
Parshin conjecture (see \cite{Gei98}).
\end{enumerate}

\subsubsection*{Counterexample for non-commutative discrete valuation rings 
(Due to Kazuya Kato)}

In Gersten's conjecture, the assumption of commutativity is essential. 
Let $D$ be a skew field finite over $\bQ_p$, 
$A$ its integer ring and $a$ its prime element. 
As the inner automorphism of $a$ induces non-trivial automorphism on 
its residue field, we have $x\in A^{\times}$ with $y=axa^{-1}x^{-1}$ 
is non vanishing in its residue field, a fortiori in 
$K_1(A)={(A^{\times})}^{\ab}$. 
On the other hand, $y$ is a commutator in $D^{\times}$. 
Hence it turns out that the canonical map 
$K_1(A)\to K_1(D)={(D^{\times})}^{\ab}$ is not injective. 

\subsection*{Idea of the proof}

My idea of how to prove Gersten's conjecture has come from weight argument 
of Adams operations in \cite{GS87} and \cite{GS99}. 
In my viewpoint, difficulty of solving Gersten's conjecture consists of 
ring theoretic side and homotopy theoretic side. 
We will explain ideas about how to overcome each difficulty.

\subsubsection*{Ring theoretic side}
\label{subsubsec:Ring theoretic side}

Combining with the results in \cite{GS87} and \cite{TT90}, 
for a commutative regular local ring $R$, there exists 
Adams operations $\{\varphi_k\}_{k\geq 0}$ on 
$K_0(\cM_R^p)$ and we have the equality
\begin{equation}
\label{eq:Adams operation}
\varphi_k([R/(f_1,\cdots,f_p)])=k^p[R/(f_1,\cdots,f_p)]
\end{equation}
where a sequence $f_1,\cdots,f_p$ is an $R$-regular sequence 
(see \cite[Proposition 2.2]{Moc17a}). 
Thus roughly saying, the generator conjecture says 
that for each $p$, $K_0(\cM_R^p)$ is spanned by 
objects of weight $p$ and Gersten's conjecture 
could follow from weight argument of Adams operations. 
We illustrate how to prove that the generator conjecture implies 
Gersten's conjecture for $K_0$ without using Adams operations. 
\begin{proof}
Let a sequence $f_1,\cdots,f_p$ be an $R$-regular sequence. 
Then there exists the short exact sequence
$$0\to R/(f_1,\cdots,f_{p-1})\onto{f_p} R/(f_1,\cdots,f_{p-1}) \to R/(f_1,\cdots,f_{p})\to 0$$
in $\cM_R^{p-1}$. 
Thus the class $ $ in $K_0(\cM_R^p)$ goes to 
$$[R/(f_1,\cdots,f_{p})]=[R/(f_1,\cdots,f_{p-1})]-[R/(f_1,\cdots,f_{p-1})]=0$$
in $K_0(\cM_R^{p-1})$.
\end{proof}

\paragraph{Gambit A}
We will establish and prove a higher analogue of the generator conjecture.

\sn
Inspired from the works 
\cite{Iwa59}, \cite{Ser59}, \cite{Bou64}, \cite{Die86} and \cite{Gra92}, 
we establish a classification theory of modules by utilizing cubes 
in \cite{Moc13a}, \cite{Moc13b} and \cite{MY14}. 
In this article, we will implicitly use these theory 
and simplify the arguments in \cite{Moc13a} and \cite{Moc15}.

\subsubsection*{Homotopy theoretic side}
\label{subsubsec:Homotopy theoretic side}

Roughly saying, we will trying to compare the following two functors 
on $K$-theory. 
We denote the category of bounded chain complexes on $\cM_R^{p-1}$ 
by $\Ch_b(\cM_R^{p-1})$. 
$$\cM_R^p\to \Ch_b(\cM_R^{p-1}),$$
$$R/(f_1,\cdots,f_p)\mapsto
\begin{cases}
\begin{bmatrix}R/(f_1,\cdots,f_{p-1})\\ \downarrow f_p\\ R/(f_1,\cdots,f_{p-1}) \end{bmatrix} \underset{\qis}{\sim} R/(f_1,\cdots,f_p)\\
\begin{bmatrix} R/(f_1,\cdots,f_{p-1})\\ \downarrow \id\\ R/(f_1,\cdots,f_{p-1})\end{bmatrix} \underset{\qis}{\sim} 0.
\end{cases}
 $$
The functors above shall be homotopic to each other on $K$-theory 
by the additivity theorem. 
A problem is that the functors above are not {\it $1$-functorial}!! 
We need to a suitable notion of $K$-theory for higher categories 
or need to discuss more subtle argument for such exotic functors. 

\paragraph{Gambit B}
We give a modified definition of algebraic $K$-theory in a particular 
situation and establish a technique of rectifying lax functors 
to $1$-functors and by utilizing this definition and these techniques, 
we will treat such exotic functors inside the classical Waldhausen $K$-theory.

\subsection*{Local Gersten's conjecture}
\label{subsec:Local Gersten's conjecture}

Let $A$ be a commutative noetherian ring with $1$ and 
let $p$ be a positive integer and let $I$ 
be an ideal of $A$ with codimension $V(I)\geq p$ in $\Spec A$. 
Let $\cM_A^I$ be a full subcategory of $\cM_A^p$ consisting of those 
modules $M$ supported on $V(I)$ and 
$\cM^I_{A,\red}$ be a full subcategory of $\cM_A^I$ consisting of 
those modules $M$ such that $IM$ are trivial. 
For a non-negative integer $q$, 
we denote the full subcategory of $\cM_A^I$ and $\cM_{A,\red}^I$ consisting of those $A$-modules $M$ with $\pd_AM\leq q$ by $\cM_A^I(q)$ and 
$\cM_{A,\red}^I(q)$ respectively. 
Let $\cP_A$ be the category of finitely generated projective $A$-modules. 

\sn
For a commutative regular local ring $R$, let $J$ be an ideal generated 
by $R$-regular sequence $f_1,\cdots,f_p$ such that $f_i$ is a prime 
element for any $1\leq i\leq p$. 
First notice the following isomorphisms for each non-negative integer $n$:
\begin{equation}
\label{eq:colimit description}
K_n(\cM_R^p)\underset{\textbf{I}}{\isoto}\underset{\substack{\codim_{\Spec R} V(I)=p \\ \Spec R/I\underset{\text{regular}}{\rinc} \Spec R }}{\colim}K_n(\cM_R^I),\end{equation}
\begin{equation}
\label{eq:resolution}
K_n(\cM^J_{R,\red})\underset{\textbf{II}}{\simeq}K_n(\cM_R^J).
\end{equation}
Since $R$ is Cohen-Macaulay, the 
ordered set of all ideals of $R$ that contains an $R$-regular sequence of 
length $p$ with usual inclusion is directed. 
Thus $\cM_R^p$ is the filtered limit $\displaystyle{\colim_I\cM_R^I}$ 
where $I$ runs through any ideal generated by 
any $R$-regular sequence of length $p$. 
Thus the isomorphism $\textbf{I}$ follows from co-continuity of 
$K$-theory. 
The isomorphism $\textbf{II}$ follows from the d\'evissage theorem. 
Hence we propose the following local Gersten's conjecture. 
Let $\ff_S=\{f_s\}_{s\in S}$ be a family of elements 
indexed by a non-empty finite set $S$ in $A$. 
For a non-empty subset $T\subset S$, 
we denote $\ff_T:=\{f_t\}_{t\in T}$ and $\ff_TA$ stands for an ideal of $A$ 
spanned by $\ff_T$. 
By convention, 
we set $\ff_{\emptyset}A=(0)$ the zero ideal of $A$. 
For a finite set $T$, $\#T$ stands for the cardinality of $T$.

\begin{con}[\bf Local Gersten's conjecture]
\label{con:Local Gersten's conjecture}
Let $R$ be a commutative regular local ring and let 
$\ff_S=\{f_s\}_{s\in S}$ be an $R$-regular sequence. 
For an element $s\in S$, 
the inclusion functor $\cM_{R,\red}^{\ff_SR}(\#S)\rinc 
\cM_{R}^{\ff_{S\ssm\{s\}}R}(\#S-1)$ induces the zero map on $K$-theory.
\end{con}

\sn
Local Gersten's conjecture for any regular sequences of 
any regular local rings imply original Gersten's conjecture. 
Our main theorem is the following:

\begin{thm}
\label{main theorem}
Let $B$ be a commutative noetherian local ring with $1$ and let 
$g$ be a non-zero divisor of $B$. 
Then the inclusion functor $\cP_{B/gB}\rinc \cM_B(1)$ 
induces the zero map on $K$-theory.
\end{thm}

\sn
We obtain the following result from Theorem~\ref{main theorem}.

\begin{cor}
\label{cor:LGC for rsp}
Let $R$ be a commutative regular local ring and let $\ff_S=\{f_s\}_{s\in S}$ 
be a part of regular system of parameter. 
Then local Gersten's conjecture for $R$ and $\ff_S$ is true. 
Namely the inclusion functor $\cM^{\ff_SR}_{R,\red}(\#S)\rinc \cM_{R,\red}^{\ff_{S\ssm\{s\}}R}(\#S-1)$ induces the zero map on $K$-theory.
\end{cor}

\begin{proof}
First we prove the equality 
\begin{equation}
\label{eq:regular system case cMRred}
\cM_{R,\red}^{\ff_SR}(\#S)=\cP_{R/\ff_SR}
\end{equation}
by induction on the cardinality of $S$. 
If $S=\emptyset$, assertion is trivial. 
For $\#S\geq 1$, we fix an element $s\in S$ and let $M$ be an $R$-module 
in $\cM^{\ff_SR}_{R,\red}(\#S)$. 
Since $R/f_sR$ is regular, 
$\pd_{R/f_sR}M<\infty$ and since $f_sM=0$, we have 
the equality $\depth_{R/f_sR}M=\depth_RM$. 
Therefore we have the equalities:
$$
\pd_{R/f_sR}M=\dim R/f_sR -\depth_{R/f_sR}M=
(\dim R-1)-\depth_RM=\pd_RM-1.
$$
This equalities show that $M$ is in 
$\cM_{R/f_sR,\red}^{\ff_{S\ssm\{s \}R}}(\#S-1)$ and this category is 
equal to $\cP_{R/\ff_SR}$ by inductive hypothesis. 
Hence we obtain the equality $\mathrm{(\ref{eq:regular system case cMRred})}$. 

\sn
Next notice that the inclusion functor 
$\cP_{R/\ff_SR}\rinc \cM_{R,\red}^{\ff_{S\ssm\{s\}}R}(\#S-1)$ 
factors through 
$$\cP_{R/\ff_SR}\rinc \cM_{R/\ff_{S\ssm\{s\}}R}(1) \rinc \cM_{R,\red}^{\ff_{S\ssm\{s\}}R}(\#S-1).$$
Applying Theorem~\ref{main theorem} to $B=R/\ff_{S\ssm\{s\}}R$ and $g=f_s$, 
we obtain the result. 
\end{proof}

\begin{proof}[Proof of Corollary~\ref{main theorem cor 1}]
Let $\ff_S=\{f_s\}_{s\in S}$ be a regular system of parameter of $R$. 
Then we have $\cM_R^d=\cM_R^{\ff_SR}(\#S)$. 
Hence we obtain the result from Corollary~\ref{cor:LGC for rsp}. 
\end{proof}

\subsection*{Strategy of the proof for main theorem}
\label{subsec:strategy of the proof}

We let $\Ch_b(\cM_B(1))$ denote the category of 
bounded complexes on $\cM_B(1)$ 
(we use homological notation). 
Let $\cC$ be the full subcategory of $\Ch_b(\cM_B(1))$ consisting of those 
complexes $x$ such that $x_i=0$ unless $i=0$, $1$ and $x_1$ and $x_0$ are 
free $B$-modules and the bounded map $d^x\colon x_1\to x_0$ is injective 
and $\Homo_0x:=\coker(x_1\onto{d^x}x_0)$ is annihilated by $g$. 

\begin{lem}
\label{lem:exact structure of C}
The category $\cC$ naturally becomes an idempotent complete 
exact category such that 
the inclusion functor $\eta\colon\cC\rinc \Ch_b(\cM_B(1))$ is exact and 
reflects exactness. 
Moreover the functor $\Homo_0\colon\cC\to \cP_{B/gB}$ is exact. 
\end{lem}

\begin{proof}
Let $\cC'$ be the full subcategory of $\Ch_b(\cM_B(1))$ consisting of those 
complexes $x$ such that $x_i=0$ unless $i=0$, $1$ and $x_1$ and $x_0$ 
are free $B$-modules and the bounded map 
$d^x\colon x_1\to x_0$ is injective and 
$\Homo_0x:=\coker(x_1\onto{d^x}x_0)$ is annihilated by $g^n$ for 
some positive integer $n$. 
Then $\cC'$ is closed under extensions in $\Ch_b(\cM_B(1))$. 
Thus $\cC'$ is naturally becomes an exact category such that the inclusion 
functor $\cC'\rinc \Ch_b(\cM_B(1))$ is exact and reflects exactness. 
Moreover $\cC$ is closed under admissible sub- and quotient objects and 
finite direct sums  and direct summands in $\cC'$. 
Hence $\cC$ is naturally becomes an exact category such that the inclusion 
functor $\cC\rinc \cC'$ is exact and reflects exactness. 
Since boundary maps of objects in $\cC$ are injective, 
the functor $\Homo_0$ from $\cC$ to the category of 
finitely generated projective $B/gB$-modules is exact. 
\end{proof}

\sn
Thus we obtain the commutative diagram
$$
\xymatrix{
K(\cC) \ar[r]^{\!\!\!\!\!\!\!\!\!\!\!\!\!\!\!\!\!\!\!\!\!\!\!\!\!\!K(\eta)} \ar[d]_{K(\Homo_0)} & K(\Ch_b(\cM_B(1));\qis)\\
K(\cP_{B/gB}) \ar[r] & K(\cM_B(1)) \ar[u]^{\wr}_{\textbf{I}}
}
$$
where $\qis$ is the class of all quasi-isomorphisms in $\Ch_b(\cM_B(1))$ 
and the map $\textbf{I}$ which is induced from the inclusion functor 
$\cM_B(1)\rinc \Ch_b(\cM_B(1))$ is a homotopy equivalence by 
Gillet-Waldhausen theorem.

We will prove that 
\begin{enumerate}
\item[$\mathrm{(\alpha)}$]
The map $K(\Homo_0)$ is a split epimorphism in the stable category of spectra 
(see \S~\ref{sec:A proof of assertion alpha}) and

\item[$\mathrm{(\beta)}$]
The map $K(\eta)$ is the zero map in the stable category of spectra 
(see \S~\ref{sec:A proof of assertion beta}).
\end{enumerate}

\sn
Then we will complete the proof of Theorem~\ref{main theorem}. 
Assertion $\mathrm{(\alpha)}$ corresponds with Gambit A and 
assertion $\mathrm{(\beta)}$ corresponds with Gambit B in the previous 
subsection. 

\section{A proof of assertion $\mathrm{(\alpha)}$}
\label{sec:A proof of assertion alpha}

In this section, we will give a proof of assertion that $K(\Homo_0)$ is 
a split epimorphism in the stable category of spectra. 

\sn
Let $\calD$ be the full subcategory of $\Ch_b(\cM_B(1))$ consisting of those 
complexes $x$ such that $x_i=0$ unless $i=0$, $1$ and the bounded 
map $d^x\colon x_1\to x_0$ is injective and 
$\Homo_0x:=\coker(x_1\onto{d^x}x_0)$ is in $\cM_{B,\red}^{gB}(1)$ 
where $gB$ is an ideal of $B$ spanned by $g$. 
As similar to a proof for $\cC$ 
(Lemma~\ref{lem:exact structure of C}), 
we can show that $\calD$ is an exact category such that the inclusion functor 
$\calD\rinc \Ch_b(\cM_B(1))$ is exact and reflects exactness and we can also 
show that the functor $\Homo_0\colon\calD\to \cM_{B,\red}^{gB}(1)$ is 
exact. 
Thus we obtain the commutative square below
$$
\xymatrix{
K(\cC) \ar[r] \ar[d]_{K(\Homo_0)} & K(\calD) \ar[d]^{K(\Homo_0)}\\
K(\cP_{B/gB}) \ar[r] & K(\cM_{B,\red}^{gB}(1)).
}
$$
where the horizontal maps are induced from the inclusion functors. 
Since the functor $\Homo_0\colon\calD\to\cM_{B,\red}^{gB}(1)$ 
admits a section which is defined by sending an object $x$ 
in $\cM_{B,\red}^{gB}(1)$ to the complex $[0\to x]$ in $\calD$, 
the right vertical map in the diagram above is a split epimorphism. 

\begin{lem}
\label{lem:resolution theorem}
{\rm (}Compare {\rm \cite[2.21]{Moc13a}} and {\rm\cite[2.1.1]{Moc17b}).}\ \ 
The inclusion functors $\cC\rinc\calD$ and $\cP_{B/gB}\rinc\cM_{B,\red}^{gB}(1)$ induce homotopy equivalences of spectra on $K$-theory respectively. 
\end{lem}

\begin{proof}
We apply Quillen's resolution theorem to inclusion functors 
$\cC\rinc\calD$ and $\cP_{B/gB}\rinc\cM_{B,\red}^{gB}(1)$ respectively. 
Thus what we need to check are the following three conditions:
\begin{enumerate}
\item
$\cP_{B/gB}$ (resp. $\cC$) is closed under extensions 
in $\cM_{B,\red}^{gB}(1)$ (resp. $y\in\Ob\calD$).

\item
For an admissible exact sequence
\begin{equation}
\label{eq:admissible exact sequence}
z\rinf y\rdef x,
\end{equation}
in $\cM_{B,\red}^{gB}(1)$ (resp. $\calD$), if $y$ is in $\cP_{B/gB}$ 
(resp. $\cC$), then $z$ is also in $\cP_{B/gB}$ (resp. $\cC$). 
\item
For an object $x$ in $\cM_{B,\red}^{gB}(1)$ (resp. $\calD$), 
there exists an admissible exact sequence $\mathrm{(\ref{eq:admissible exact sequence})}$ in $\cM_{B,\red}^{gB}(1)$ (resp. $\calD$) with 
$y\in\Ob\cP_{B/gB}$ (resp. $y\in\Ob\cC$).
\end{enumerate}
First 
notice that since $B$ is local, 
every finitely generated projective $B$-modules and $B/gB$-modules are free. 
Assertion $1$ follows from semi-simplicity of $\cP_{B/gB}$ and $\cC$ 
(for semi-simplicity of $\cC$, 
see Lemma~\ref{lem:structure of C} $\mathrm{(II)}$). 

Next we show assertion $2$. 
In the admissible exact sequence 
$\mathrm{(\ref{eq:admissible exact sequence})}$, 
$B$-projective dimension(s) of $B$-module $x$ is 
(resp. $x_i$ ($i=0$, $1$ are) less than $1$, 
$z$ is a (resp. $z_i$ ($i=0$, $1$) are) finitely generated free $B/gB$-module 
(resp. $B$-modules). 

Finally we give a proof of assertion $3$ for the pair $\cC\rinc\calD$. 
For an object $x$ in $\calD$, there exists a surjection 
$y_0\rdef x_0$ with $y_0\in\Ob\cP_{B/gB}$. 
We denote the kernel of compositions $y_0\rdef x_0\rdef \Homo_0x$ by $y_1$ 
and write $d^y$ for the inclusion map $y_1\to y_0$. 
Then since $\pd_B\Homo_0x\leq 1$, $y_1$ is a projective $B$-module, 
namely a free $B$-module. 
The snake lemma shows the induced morphism $y_1\to x_1$ which 
makes diagram below commutative is surjective. 
$$
\xymatrix{
y_1 \ar@{>->}[d]_{d^y} \ar@{->>}[r] & x_1 \ar@{>->}[d]^{d^x}\\ 
y_0 \ar@{->>}[d] \ar@{->>}[r] & x_0 \ar@{->>}[d]\\
\Homo_0x \ar@{=}[r] & \Homo_0x.
}
$$
Thus $y=[y_1\onto{d^y}y_0]$ is in $\cC$. 
We denote the kernel of $y\to x$ by $z$. 
Then by the snake lemma shows $\Homo_iz=0$ for $i=0$, $1$. 
Thus $z$ is in $\calD$ and it turns out that 
the map $y\rdef x$ is 
an admissible epimorphism.
\end{proof}

Thus the horizontal maps in the diagram above are homotopy equivalences 
of spectra and the left vertical map is also a split epimorphism in 
the stable category of spectra. 

\section{A proof of assertion $\mathrm{(\beta)}$}
\label{sec:A proof of assertion beta}

In this section we will give a proof of assertion that 
$K(\eta)$ is the zero map in the stable category of spectra. 
Let $\cB$ be the full subcategory of $\Ch_b(\cM_B(1))$ consisting of those 
complexes $x$ such that $x_i=0$ unless $i=$ or $i=1$. 
Let $\fs_i\colon\cB\to\cM_B(1)$ ($i=0$, $1$) be an exact functor 
defined by sending an object $x$ in $\cB$ to $x_i$ in $\cM_B(1)$. 
By the additivity theorem, the map $\fs_1\times\fs_2\colon iS_{\cdot}\cB\to iS_{\cdot}\cM_B(1)$ is a homotopy equivalence. 
Let $j\colon \cB\to \Ch_b(\cM_B(1))$ be the inclusion functor. 
We wish to define two exact `functor' $\mu_1$, $\mu_2\colon\cC\to\cB$ which 
satisfy the following conditions:
\begin{enumerate}
\item[$\mathrm{(A)}$]
We have the equality
\begin{equation}
\label{eq:s1timess2}
\fs_1\times\fs_2\mu_1=\fs_1\times\fs_2\mu_2.
\end{equation}
\item[$\mathrm{(B)}$]
The are natural transformations $\eta\to j\mu_1$ and $0\to j\mu_2$ such that 
all components are quasi-isomorphisms. 
\end{enumerate} 

Then we have the equalities
$$ 
K(\eta)=K(j\mu_1)=K(j){K(\fs_1\times\fs_2)}^{-1}K(\fs_1\times\fs_2\mu_1)=
K(j){K(\fs_1\times\fs_2)}^{-1}K(\fs_1\times\fs_2\mu_2)=K(j\mu_2)=0.
$$
To define the `functors' $\mu_i$ ($i=1$, $2$), we analyze morphisms in $\cC$.

\subsection*{Structure of $\cC$}
\label{subsec:structure of C}

\begin{df}
\label{df:(n,m)_B} (Compare \cite[1.1.3]{Moc17b}).\ \ 
For a pair of non-negative integers $n$ and $m$, we write 
${(n,m)}_B$ for the complex of the form 
$\begin{bmatrix}
\!\!\!\!\!\!\!\!\!\!\!\!\!\!\!\!\!\!\!\!\!\!\!\!\!\!B^{\oplus n}\oplus B^{\oplus m}\\ 
\ \ \ \ \ \ \ \ \ \downarrow \begin{pmatrix}gE_n & 0 \\ 0 & E_m \end{pmatrix}\\ 
\!\!\!\!\!\!\!\!\!\!\!\!\!\!\!\!\!\!\!\!\!\!\!\!\!\!B^{\oplus n}\oplus B^{\oplus m}
\end{bmatrix}$
in $\cC$ where $E_k$ is the $k\times k$ unit matrix. 
\end{df}

\begin{lem}
\label{lem:structure of C}
\begin{enumerate}
\enumidef

\item
{\rm (Compare \cite[1.2.10, 1.2.13]{Moc17b}).}\ \ 
Let $n$ be a positive integer. 
For any endomorphism $a\colon {(n,0)}_B \to {(n,0)}_B$, 
the following conditions are equivalent.
\begin{enumerate}
\enumiidef
\item
$a$ is an isomorphism.
\item
$a$ is a quasi-isomorphism.
\end{enumerate}

\item
{\rm (Compare \cite[2.17]{Moc13a}).}\ \ 
An object in $\cC$ is projective. 
In particular, $\cC$ is a semi-simple exact category.

\item
{\rm (Compare \cite[1.2.15]{Moc17b}).}\ \ 
For any object $x$ in $\cC$, 
there exists a pair of non-negative integers $n$ and $m$ 
such that $x$ is isomorphic to ${(n,m)}_B$.
\end{enumerate}
\end{lem}

\begin{proof}
$\mathrm{(I)}$ 
We assume condition $\mathrm{(b)}$. 
In the commutative diagram below
$$
\xymatrix{
B^{\oplus n} \ar[r]^{gE_n} \ar[d]_{a_1} & B^{\oplus n} \ar[r] \ar[d]_{a_0} 
& \Homo_0({(n,0)}_B) \ar[d]^{\Homo_0(a)}\\
B^{\oplus n} \ar[r]_{gE_n} & B^{\oplus n} \ar[r] & 
\Homo_0({(n,0)}_B), 
}
$$
first we will prove that $a_0$ is an isomorphism. 
Then $a_1$ is also an isomorphism by the five lemma. 
By taking determinant 
of $a_0$, we shall assume that $n=1$. 
Then assertion follows from Nakayama's lemma.

\sn
$\mathrm{(II)}$ 
Let $t\colon y\rdef z$ be an admissible epimorphism in $\cC$ 
and let $f\colon x\to z$ be a morphism in $\cC$. 
Then since $\Homo_0(x)$ is a projective $B/gB$-module, 
there exists a homomorphism of $B/gB$-modules 
$\sigma\colon \Homo_0(x)\to\Homo_0(y)$ 
such that $\Homo_0(t)\sigma =\Homo_0(f)$. 
Since $x$ is a complex of free $B$-modules, there exists a 
morphism of complexes $s'\colon x\to y$ such that $\Homo_0(s')=\sigma$ 
and $ts'$ is chain homotopic to $f$ by 
\cite[Comparison theorem 2.2.6]{Wei94}. 
Namely there exists a map $h\colon x_0\to z_1$ such that 
${(f-ts')}_0=d^zh$ and ${(f-ts')}_1=hd^x$.
$$
\xymatrix{
x_1 \ar[r]^{{(f-ts')}_1} \ar[d]_{d^x} & z_1 \ar[d]^{d^z}\\
x_0 \ar[ur]_h \ar[r]_{{(f-ts')}_0} & z_0.
}
$$
Since $x_0$ is projective, there exists a map $u\colon x_0\to y_1$ 
such that $t_1u=h$. 
We set $s_1:=s'_1+ud^x$ and $s_0:=s'_0+d^yu$. 
Then we can check that $s$ is a morphism of complexes of 
$B$-modules and $f=ts$. 

\sn
$\mathrm{(III)}$ 
By considering $\displaystyle{x\otimes_B B\left [\frac{1}{g}\right ]}$, 
it turns out that $x_1$ and $x_0$ are same rank. 
Thus we shall assume $x_1=x_0=B^{\oplus m}$. 
First assume that $x$ is acyclic. 
Then the boundary map $d^x\colon x_1\to x_0$ is invertible and 
$$ 
\begin{bmatrix}
x_1\\ \ \ \ \downarrow d^x \\ x_0 
\end{bmatrix}
\begin{matrix}
\onto{d^x}\\ \ \ \\ \overset{\to}{\id}
\end{matrix}
\begin{bmatrix}
B^{\oplus m}\\ \ \ \ \ \ \downarrow \id_{B^{\oplus m}}\\ B^{\oplus m}
\end{bmatrix}
$$
gives an isomorphism between $x$ and ${(0,m)}_B$. 
Thus we obtain the result in this case. 

Next assume that $\Homo_0(x)\neq 0$. 
Then since $\Homo_0(x)$ is a finitely generated $B/gB$-modules, 
there exists a positive integer $n$ and an isomorphism 
$\sigma \colon {(B/gB)}^{\oplus n} \isoto \Homo_0(x)$. 
Then by \cite[Comparison theorem 2.2.6]{Wei94}, 
there exists morphisms of complexes in $\cC$, 
${(n,0)}_B\onto{a} x$ and $x\onto{b} {(n,0)}_B$ such that $\Homo_0(a)=\sigma$ and $\Homo_0(b)=\sigma^{-1}$. 
Thus by $\mathrm{(I)}$, $ba$ is an isomorphism. 
By replacing $a$ with $a{(ba)}^{-1}$, 
we shall assume that $ba=\id$. 
Hence there exists a complex $y$ in $\cC$ and a split exact sequence:
\begin{equation}
\label{eq:split exact seq}
{(n,0)}_B\overset{a}{\rinf} x\rdef y. 
\end{equation}
Since boundary maps of objects in $\cC$ are injective, 
the functor $\Homo_0$ from $\cC$ to the category of 
finitely generated projective $B/gB$-modules is exact. 
By taking $\Homo_0$ to the sequence 
$\mathrm{(\ref{eq:split exact seq})}$, 
it turns out that $y$ is acyclic and by the first paragraph, 
we shall assume  that $y$ is of the form ${(0,m)}_B$. 
Hence $x$ is isomorphic to ${(n,m)}_B$.
\end{proof}

\begin{conv}
\label{conv:matrix representations}
(Compare \cite[1.2.16]{Moc17b}).\ \ 
We can denote a morphism $\varphi\colon{(n,m)}_B\to{(n',m')}_B$ of $\cC$ by
$$
\begin{bmatrix}
\!\!\!\!\!\!\!\!\!\!\!\!\!\!\!\!\!\!B^{\oplus n}\oplus B^{\oplus m}\\ 
\ \ \ \ \ \ \ \ \ \ \ \ \ \downarrow \begin{pmatrix}gE_n & 0 \\ 0 & E_m \end{pmatrix}\\
\!\!\!\!\!\!\!\!\!\!\!\!\!\!\!\!\!\!B^{\oplus n}\oplus B^{\oplus m}
\end{bmatrix}
\begin{matrix}
\overset{\varphi_1}{\to}\\ \ \ \\ \underset{\varphi_0}{\to}
\end{matrix}
\begin{bmatrix}
\!\!\!\!\!\!\!\!\!\!\!\!\!\!\!\!\!\!\!\!\!\!\!B^{\oplus n'}\oplus B^{\oplus m'}\\ 
\ \ \ \ \ \ \ \ \ \ \ \ \ \downarrow \begin{pmatrix}gE_{n'} & 0 \\ 0 & E_{m'} \end{pmatrix}\\
\!\!\!\!\!\!\!\!\!\!\!\!\!\!\!\!\!\!\!\!\!\!\!B^{\oplus n'}\oplus B^{\oplus m'}
\end{bmatrix}
$$
with 
$\displaystyle{\varphi_1=\begin{pmatrix}\varphi_{(n',n)} & \varphi_{(n',m)}\\
g\varphi_{(m',n)} & \varphi_{(m',m)} \end{pmatrix}}$ and 
$\displaystyle{\varphi_0=\begin{pmatrix}
\varphi_{(n',n)} & g\varphi_{(n',m)}\\
\varphi_{(m',n)} & \varphi_{(m',m)}\end{pmatrix}}$ 
where $\varphi_{(i,j)}$ are $i\times j$ 
matrices whose coefficients are in $B$. 
In this case we write 
\begin{equation}
\label{eq:representation of varphi}
\begin{pmatrix}
\varphi_{(n',n)} & \varphi_{(n',m)}\\
\varphi_{(m',n)} & \varphi_{(m',m)}
\end{pmatrix}
\end{equation}
for $\varphi$. 
In this matrix presentation of morphisms, the composition of 
morphisms between objects ${(n,m)}_B\onto{\varphi}{(n',m')}_B\onto{\psi}{(n'',m'')}_B$ in $\cC$ is described by 
\begin{equation}
\label{eq:composition of morphisms}
{\footnotesize{
\begin{pmatrix}
\psi_{(n'',n')} & \psi_{(n'',m')} \\ 
\psi_{(m'',n')} &\psi_{(m'',m')}
\end{pmatrix}
\begin{pmatrix} 
\varphi_{(n',n)} & \varphi_{(n',m)}\\
\varphi_{(m',n)} & \varphi_{(m',m)}
\end{pmatrix}=
\begin{pmatrix}
\psi_{(n'',n')}\varphi_{(n',n)}+g\psi_{(n'',m')}\varphi_{(m',n)} & 
\psi_{(n'',n')}\varphi_{(n',m)}+\psi_{(n'',m')}\varphi_{(m',m)}\\
\psi_{(m'',n')}\varphi_{(n',n)}+\psi_{(m'',m')}\varphi_{(m',n)} & 
g\psi_{(m'',n')}\varphi_{(n',m)}+\psi_{(m'',m')}\varphi_{(m',m)}\\
\end{pmatrix}.
}}
\end{equation}
Thus the category $\cC$ is categorical equivalent to the category 
whose objects are ordered pair of non-negative integers $(n,m)$ 
and whose morphisms from an object $(n,m)$ to 
$(n',m')$ are $2\times 2$ matrices 
of the form $\mathrm{(\ref{eq:representation of varphi})}$ 
of $i\times j$ matrices 
$\varphi_{(i,j)}$ whose coefficients are in $B$ 
and compositions are given by the formula 
$\mathrm{(\ref{eq:composition of morphisms})}$. 
We sometimes identify these two categories.
\end{conv}

\subsection*{Modified algebraic $K$-theory}
\label{subsec:Modified algebraic K-theory}

A candidate of a pair of $\mu_i\colon\cC\to \cB$ ($i=1$, $2$) are 
following. 
For an element $h$ in $B$, we write $\Typ(h)$ for the 
complex $[B\onto{h}B]$ in $\Ch_{[0,1]}(\cM_B)$ the category of those 
complexes $x$ such that $x_i=0$ unless $i=0$, $1$. 
We define $\mu_1'$, $\mu_2'\colon \cC\to \cB$ to be association 
by sending an object ${(n,m)}_B$ in $\cC$ to 
${\Typ(g)}^{\oplus n}$ and $\Typ(1)^{\oplus n}$ respectively and 
a morphism\\
$\displaystyle{\varphi=  
\begin{pmatrix}
\varphi_{(n',n)} & \varphi_{(n',m)}\\
\varphi_{(m',n)} & \varphi_{(m',m)}
\end{pmatrix}
\colon {(n,m)}_B\to {(n',m')}_B}$ in $\cC$ to 
$\begin{bmatrix} 
B^{\oplus n}\\
\downarrow gE_n\\
B^{\oplus n}
\end{bmatrix}
\begin{matrix}
\overset{\varphi_{(n',n)}}{\to}\\
\ \\
\underset{\varphi_{(n',n)}}{\to}
\end{matrix}
\begin{bmatrix}
B^{\oplus n}\\
\downarrow gE_{n'}\\
B^{\oplus n}
\end{bmatrix}$ and 
$
\begin{bmatrix} 
B^{\oplus n}\\
\downarrow E_n\\
B^{\oplus n}
\end{bmatrix}
\begin{matrix}
\overset{\varphi_{(n',n)}}{\to}\\
\ \\
\underset{\varphi_{(n',n)}}{\to}
\end{matrix}
\begin{bmatrix}
B^{\oplus n}\\
\downarrow E_{n'}\\
B^{\oplus n}
\end{bmatrix}
$ respectively. 
Notice that they are not $1$-functors. 
We need to make revision in the previous idea. 
We introduce a modified version of algebraic $K$-theory of $\cC$. 

\begin{df}[\bf Triangular morphisms]
\label{df:triangular morphisms}
(Compare \cite[2.2.5]{Moc17b}).\ \ 
We say that a morphism $\varphi\colon {(n,m)}_B\to {(n',m')}_B$ 
in $\cC$ of the form 
$\mathrm{(\ref{eq:representation of varphi})}$ 
in an {\it upper triangular} if $\varphi_{(m',n)}$ 
is the zero morphism, and say that $\varphi$ is a 
{\it lower triangular} if $\varphi_{(n',m)}$ is the zero morphism. 
We denote the class of all lower triangular isomorphisms in $\cC$ by 
$i^{\bigtriangledown}$. 
Next we define $S_{\cdot}^{\bigtriangleup}\cC$ to be a simplicial subcategory 
of $S_{\cdot}\cC$ consisting of those objects $x$ such that 
$x(i\leq j)\to x(i'\leq j')$ is an upper triangular morphism for 
each $i\leq i',j\leq j'$. 
\end{df}

\begin{lem}
\label{lem:modified K-theory}
\begin{enumerate}
\enumidef
\item
{\rm(}Compare {\rm \cite[2.3.5]{Moc17b}).}\ \ 
The inclusion functor 
$k\colon iS_{\cdot}^{\bigtriangleup}\cC\to iS_{\cdot}\cC$ 
is a homotopy equivalence. 
\item
{\rm(}Compare {\rm \cite[2.3.6]{Moc17b}).}\ \ 
The inclusion map 
$i^{\bigtriangledown}S_{\cdot}^{\bigtriangleup}\cC\to iS_{\cdot}\cC$ 
is a split epimorphism up to homotopy. 

\item
{\rm(}Compare {\rm \cite[2.3.7]{Moc17b}).}\ \ 
The associations $\mu_1'$ and $\mu_2'$ induce the simplicial functors 
$\mu_1$, $\mu_2\colon i^{\bigtriangledown}S^{\bigtriangleup}_{\cdot}\cC
\to iS_{\cdot}\cB$. 

\item
{\rm(}Compare {\rm \cite[2.3.7]{Moc17b}).}\ \ 
$\mu_1$ is homotopic to $\mu_2$.
\end{enumerate}
\end{lem}

\begin{proof}
$\mathrm{(I)}$ 
Since $\cC$ is semi-simple by Lemma~\ref{lem:structure of C} 
$\mathrm{(II)}$, 
the inclusion functor 
$k\colon iS_{\cdot}^{\bigtriangleup}\cC\to iS_{\cdot}\cC$ 
is an equivalence of categories for each degree. 
Therefore the inclusion functor $k$ induces a weak homotopy equivalence 
$NiS^{\bigtriangleup}_{\cdot}\cC\to NiS_{\cdot}\cC$.

\sn
$\mathrm{(II)}$ 
First by \cite[\S 1.4 Corollary]{Wal85}, composition 
$s_{\cdot}\cC\to i^{\bigtriangledown}S_{\cdot}\cC\to iS_{\cdot}\cC$ is 
a homotopy equivalence. 
Thus the inclusion functor 
$i^{\bigtriangledown}S_{\cdot}\cC\to iS_{\cdot}\cC$ 
is a split epimorphism up to homotopy. 

Next we will show that for a non-negative integer $n$, 
the inclusion functor 
$i^{\bigtriangledown}S^{\bigtriangleup}_n\cC\to i^{\bigtriangledown}S_n\cC$ 
is an equivalence of categories. 
Since we have the equality 
$i^{\bigtriangledown}S^{\bigtriangleup}_n\cC= i^{\bigtriangledown}S_n\cC $ 
for $n=0$, $1$, 
we will assume $n\geq 2$. 

For a pair of integers $0\leq q\leq p\leq n$, 
let ${(i^{\bigtriangledown}S_n\cC)}_{p,q}$ be the full subcategory of 
$i^{\bigtriangledown}S_n\cC$ consisting of those objects $x$ such that 
$x(q\leq i)\rinf x(q\leq i+1)$ and $x(q\leq i)\rdef x(q+1\leq i)$ for 
$p\leq i\leq n-1$ and $x(i\leq j)\rinf x(i\leq j+1)$ for 
$q+1\leq i\leq j\leq n-1$ and $x(i\leq j)\rdef x(i+1\leq j)$ for 
$q<i<j\leq n$ are upper triangular. 
Then there exists the inclusion functors
$$i^{\bigtriangledown}S^{\bigtriangleup}_n\cC={(i^{\bigtriangledown}S_n\cC)}_{1,0}\rinc {(i^{\bigtriangledown}S_n\cC)}_{2,0}\rinc\cdots \rinc 
{(i^{\bigtriangledown}S_n\cC)}_{n,0}\rinc 
{(i^{\bigtriangledown}S_n\cC)}_{2,1}\rinc\cdots\rinc 
{(i^{\bigtriangledown}S_n\cC)}_{n,n-1}=i^{\bigtriangledown}S_n\cC.$$
We will show that the inclusion functor 
${(i^{\bigtriangledown}S_n\cC)}_{p,q}\rinc{(i^{\bigtriangledown}S_n\cC)}_{p+1,q}$ 
is essentially surjective 
for any pair of integers $0\leq q\leq p\leq n-1$. 

Let $x$ be an object in ${(i^{\bigtriangledown}S_n\cC)}_{p+1,q}$. 
We set $\alpha_x:=\UT(x(q\leq p \to q\leq p+1))$ 
and we define $x'$ to be an object in ${(i^{\bigtriangledown}S_n\cC)}_{p,q}$ 
and an isomorphism $\gamma\colon x'\isoto x$ in the following way. 
$$
x'(i\leq j)=x(i\leq j)
$$
$$
x'(i\leq j \to i\leq j+1)=
\begin{cases}
\alpha_x^{-1}x(q\leq p-1\to q\leq p) & \text{if $(i,j)=(q,p-1)$}\\
x(q\leq p \to q\leq p+1)\alpha_x & \text{if $(i,j)=(q,p)$}\\
x(i\leq j \to i\leq j+1) & \text{otherwise}
\end{cases}
$$
$$
x'(i<j \to i+1\leq j)=
\begin{cases}
\alpha_x^{-1}x(q-1\leq p \to q\leq p) & \text{if $(i,j)=(q-1,p)$}\\
x(q\leq p \to q+1 \leq p)\alpha_x & \text{if $(i,j)=(q,p)$}\\
x(i\leq j \to i+1\leq j) & \text{otherwise}
\end{cases}
$$
$$
\gamma(i\leq j)=
\begin{cases}
\alpha_x & \text{if $(i,j)=(q,p)$}\\
\id_{x(i\leq j)} & \text{otherwise}.
\end{cases}
$$
Since {\footnotesize{
$x'(q+1\leq p \to q+1\leq p+1)x'(q\leq p \to q+1\leq p)=
x'(q\leq p+1 \to q+1\leq p+1)x'(q\leq p \to q\leq p+1)$}} 
is upper triangular and the map 
$x'(q+1\leq p \to q+1\leq p+1)$ is a monomorphism, 
it turns out that the map 
$x'(q\leq p\to q+1\leq p)$ is also upper triangular. 
Thus $x'$ is in ${i^{\bigtriangledown}S_{\cdot}\cC}_{p,q}$ 
and we obtain the result.

Similarly we can show that the inclusion functor 
${i^{\bigtriangledown}S_{\cdot}\cC}_{n,q}\rinc 
{i^{\bigtriangledown}S_{\cdot}\cC}_{q+2,q+1}$ 
for any integers $0\leq q<p\leq n-1$ is an equivalence of categories. 
Thus the inclusion 
$i^{\bigtriangledown}S^{\bigtriangleup}_{\cdot}\cC\to 
i^{\bigtriangledown}S_{\cdot}\cC$ 
is a weak equivalence by the realization lemma 
\cite[Appendix A]{Seg74} or \cite[5.1]{Wal78}. 
Hence we complete the proof. 

\sn
$\mathrm{(III)}$ 
Notice that for a pair of composable 
morphisms in $\cC$, 
\begin{equation}
\label{eq:varphi psi}
{(n,m)}_B\onto{\varphi}{(n',m')}_B\onto{\psi}{(n'',m'')}_B,
\end{equation}
\begin{enumerate}
\item
if both $\varphi$ and $\psi$ 
are upper triangular or both $\varphi$ and $\psi$ 
are lower triangular, then we have the equality 
$\mu_i'(\varphi\psi)=\mu_i'(\varphi)\mu_i'(\psi)$ for $i=1$, $2$, 
\item
if the sequence $\mathrm{(\ref{eq:varphi psi})}$ is exact in $\cC$, the the sequence
$$
\mu_i'({(n,m)}_B)\onto{\mu_i'(\varphi)}\mu_i'({(n',m')}_B)\onto{\mu_i'(\psi)}
\mu_i'({(n'',m'')}_B)
$$
is exact in $\cB$ for $i=1$, $2$ 
by Lemma~\ref{lem:exact structures in C} below and
\item
if $\varphi$ is an isomorphism in $\cC$, then $\mu_i'(\varphi)$ is an isomorphism in $\cB$ for $i=0$, $1$ 
by Lemma~\ref{lem:isomorphisms in C} below. 
\end{enumerate}
Thus the associations $\mu_1'$ and $\mu_2'$ induce the simplicial functors 
$\mu_1$, $\mu_2\colon i^{\bigtriangledown}S^{\bigtriangleup}_{\cdot}\cC
\to iS_{\cdot}\cB$. 

\sn
$\mathrm{(IV)}$ 
Inspection shows the equality $\mathrm{(\ref{eq:s1timess2})}$. 
Hence $\mu_1$ is homotopic to $\mu_2$ by the additivity theorem. 
\end{proof}

\begin{df}[\bf upside-down involution]
\label{df:upside-dwon involution}
(Compare \cite[1.2.17]{Moc17b}). 
Let $x=[x_1\onto{d^x}x_0]$ be an object in $\cC$. 
Since $d^x$ is a monomorphism and $x_1$ and $x_0$ have the same rank, 
$d^x$ is invertible in $\displaystyle{B\left [\frac{1}{g}\right ]}$ and 
$g{d^x}^{-1}\colon x_0\to x_1$ is a morphism of $B$-modules. 
We define $\UD\colon \cC\to \cC$ to be a functor by sending 
an object $[x_1\onto{d^x}x_0 ]$ to $[x_0\onto{g{(d^x)}^{-1}}x_1]$ 
and a morphism $\varphi\colon x\to y$ to 
$$
\begin{bmatrix}
x_0\\ \downarrow g{(d^x)}^{-1}\\ x_1
\end{bmatrix}
\begin{matrix}
\overset{\varphi_0}{\to} \\ \ \ \\ \underset{\varphi_1}{\to}
\end{matrix}
\begin{bmatrix}
y_0\\ \downarrow g{(d^y)}^{-1}\\ y_1
\end{bmatrix}.
$$
Namely we have the equations:
$$\UD({(n,m)}_B)={(m,n)}_B,$$
$$
\UD\left (\begin{pmatrix}
\varphi_{(n',n)} & \varphi_{(n',m)}\\
\varphi_{(m',n)} & \varphi_{(m',m)} 
\end{pmatrix}\colon {(n,m)}_B\to {(n',m')}_B \right )=
\begin{pmatrix}
\varphi_{(m',m)} & \varphi_{(m',n)}\\ \varphi_{(n',m)} & \varphi_{(n',n)}
\end{pmatrix}.
$$
Obviously $\UD$ is an involution and an exact functor. 
We call $\UD$ the {\it upside-down involution}.
\end{df}

\begin{lem}
\label{lem:isomorphisms in C}
{\rm (}Compare {\rm \cite[1.2.18]{Moc17b}).}\ \ 
Let $\displaystyle{\varphi=
\begin{pmatrix}
\varphi_{(n,n)} & \varphi_{(n,m)}\\ 
\varphi_{(m,n)} & \varphi_{(m,m)}
\end{pmatrix}\colon {(n,m)}_B\to {(n,m)}_B}$ be an isomorphism in $\cC$. 
Then $\varphi_{(n,n)}$ and $\varphi_{(m,m)}$ are invertible. 
\end{lem}

\begin{proof}
For $\varphi_{(n,n)}$, assertion follows from 
Lemma~\ref{lem:structure of C} $\mathrm{(I)}$. 
For $\varphi_{(m,m)}$, we apply the same lemma to $\UD(\varphi)$.
\end{proof}

\begin{df}[\bf Upper triangulation]
(Compare \cite[2.2.5]{Moc17b}).\ \  
Let $$\displaystyle{\varphi=
\begin{pmatrix}
\varphi_{(n,n)} & \varphi_{(n,m)}\\ 
\varphi_{(m,n)} & \varphi_{(m,m)}
\end{pmatrix}\colon {(n,m)}_B\to {(n,m)}_B}$$
be an isomorphism in $\cC$. 
By Lemma~\ref{lem:isomorphisms in C}, 
$\varphi_{(m,m)}$ is an isomorphism. 
We define $\UT(\varphi)\colon{(n,m)}_B\to {(n,m)}_B$ 
to be a lower triangular isomorphism by the formula 
$\displaystyle{\UT(\varphi):=\begin{pmatrix}E_n & 0 \\ 
-\varphi_{(m,m)}^{-1}\varphi_{(m,n)} & E_m \end{pmatrix}}$. 
Then we have an equality
\begin{equation}
\label{eq:uppertriangular}
\varphi\UT(\varphi)=\begin{pmatrix}
\varphi_{(n,n)}-g\varphi_{(n,m)}\varphi_{(m,m)}^{-1}\varphi_{(m,n)} & \varphi_{(n,m)}\\ 0 & \varphi_{(m,m)} 
\end{pmatrix}
\end{equation}
We call $\UT(\varphi)$ the {\it upper triangulation of $\varphi$}. 
Notice that if $\varphi$ is upper triangular, then 
$\UT(\varphi)=\id_{{(n,m)}_B}$. 
\end{df}

\begin{lem}
\label{lem:exact structures in C}
{\rm (}Compare {\rm \cite[1.2.19]{Moc17b}).}\ \ 
Let
\begin{equation}
\label{eq:seq in C}
{(n,0)}_B\onto{\alpha}{(n',0)}_B\onto{\beta}{(n'',0)}_B
\end{equation}
be a sequence of morphisms in $\cC$ such that $\beta\alpha=0$. 
If the induced sequence of projective $B/gB$-modules
\begin{equation}
\label{eq:H of exact seq}
\Homo_0({(n,0)}_B)\onto{\Homo_0(\alpha)}\Homo_0({(n',0)}_B)\onto{\Homo_0(\beta)}\Homo_0({(n'',0)}_B)
\end{equation}
is exact, then the sequence $\mathrm{(\ref{eq:seq in C})}$ is 
also {\rm(}split{\rm)} exact.
\end{lem}

\begin{proof}
Since the sequence $\mathrm{(\ref{eq:H of exact seq})}$ is an 
exact sequence of projective $B/gB$-modules, 
it is a split exact sequence and hence 
$n'=n+n''$ and there exists a homomorphism of $B/gB$-modules 
$\overline{\gamma}\colon\Homo_0({(n'',0)}_B)\to\Homo_0({(n',0)}_B)$ 
such that $\Homo_0(\beta)\overline{\gamma}=\id_{\Homo_0({(n'',0)}_B)}$. 
Then by \cite[Comparison theorem 2.2.6]{Wei94}, 
there is a morphism of complexes of $B$-modules $ $ such that 
$\Homo_0(\gamma)=\overline{\gamma}$. 
Since $\beta\gamma$ is an isomorphism 
by Lemma~\ref{lem:structure of C} $\mathrm{(I)}$, 
by replacing $\gamma$ with $\gamma{(\beta\gamma)}^{-1}$, 
we shall assume that $\beta\gamma=\id_{{(n'',0)}_B}$. 
Therefore there is a commutative diagram
$$
\xymatrix{
{(n,0)}_B \ar[r]^{\alpha} \ar@{-->}[d]_{\delta} &
{(n',0)}_B \ar@{->>}[r]^{\beta} \ar@{=}[d] &
{(n'',0)}_B \ar@{=}[d]\\
{(n,0)}_B \ar@{>->}[r]_{\alpha'} &
{(n',0)}_B \ar@{->>}[r]_{\beta} &
{(n'',0)}_B
}
$$
such that the bottom line is exact. 
Here the dotted arrow $\delta$ is induced from the universality 
of $\Ker\beta$. 
By applying the functor $\Homo_0$ to the diagram above and 
by the five lemma, 
it turns out that $\Homo_0(\delta)$ is an isomorphism of 
projective $B/gB$-modules and hence $\delta$ is also an isomorphism 
by Lemma~\ref{lem:structure of C} $\mathrm{(I)}$. 
We complete the proof. 
\end{proof}

\subsection*{Homotopy natural transformations}
\label{subsec:homotopy natural transformations}

We denote the simplicial morphisms 
$i^{\bigtriangledown}S^{\bigtriangleup}_{\cdot}\cC \to
\qis S_{\cdot}\Ch_b(\cM_B)$ 
induced from the inclusion functor 
$\eta\colon\cC\rinc \Ch_b(\cM_B)$ by the same letter $\eta$. 
For simplicial functors
$$\eta,\ j\mu_1,\ j\mu_2,\ 0\colon 
i^{\bigtriangledown}S^{\bigtriangleup}_{\cdot}\cC \to
\qis S_{\cdot}\Ch_b(\cM_B),$$
there exists a canonical natural transformation $j\mu_2\to 0$ and 
we wish to define a map from $j\mu_1$ to $\eta$. 
A candidate of $j\mu_1\to \eta$ is the following.

\sn
For any object ${(n,m)}_B$ in $\cC$, we write 
$\delta_{{(n,m)}_B}\colon j\mu_1'({(n,m)}_B)\to \eta({(n,m)}_B)$ 
for the canonical inclusion\\
$\begin{bmatrix}
B^{\oplus n}\\
\downarrow gE_n\\
B^{\oplus n}
\end{bmatrix}
\begin{matrix}
\overset{\footnotesize{\begin{pmatrix}E_n\\ 0 \end{pmatrix}}}{\to}\\
\ \\
\underset{\footnotesize{\begin{pmatrix}E_n\\ 0 \end{pmatrix}}}{\to}
\end{matrix}
\begin{bmatrix}
\!\!\!\!\!\!\!\!\!\!\!\!\!\!\!B^{\oplus n}\oplus B^{\oplus m}\\
\ \ \ \ \ \ \ \ \ \ \ \downarrow\footnotesize{\begin{pmatrix}gE_n & 0 \\ 0 & E_m \end{pmatrix}}\\
\!\!\!\!\!\!\!\!\!\!\!\!\!\!\!B^{\oplus n}\oplus B^{\oplus m}
\end{bmatrix}$. 
This is not a simplicial natural transformation. 
But it has the following nice properties. 
\begin{enumerate}
\item
$\delta_{{(n,m)}_B}$ is a chain homotopy equivalence,

\sn
for a morphism $\displaystyle{\varphi=
\begin{pmatrix}\varphi_{(n',n)} & \varphi_{(n',m)}\\ 
\varphi_{(m',n)} & \varphi_{(m',m)}
\end{pmatrix}\colon
{(n,m)}_B\to {(n',m')}_B}$ in $\cC$,

\item
if $\varphi$ is upper triangular, we have the equality 
$\eta(\varphi)\delta_{{(n,m)}_B}=\delta_{{(n',m')}_B}j\mu_1'(\varphi)$,
\item
if $\varphi$ is lower triangular, there exists a unique homotopy 
between $\eta(\varphi)\delta_{{(n,m)}_B}$ and $\delta_{{(n',m')}_B}j\mu_1'(\varphi)$. 
Namely since we have the equality 
$\displaystyle{\eta(\varphi)\delta_{{(n,m)}_B} -\delta_{{(n',m')}_B}j\mu_1'(\varphi)=\begin{pmatrix}0 \\ \varphi_{(m',n)}\end{pmatrix}}$, the map 
$$\displaystyle{H:=\begin{pmatrix}0 \\ \varphi_{(m',n)}\end{pmatrix} \colon B^{\oplus n}\to B^{\oplus n'}\oplus B^{\oplus m'}}$$
gives a chain homotopy between $\eta(\varphi)\delta_{{(n,m)}_B}$ and $\delta_{{(n',m')}_B}j\mu_1'(\varphi)$. 
\end{enumerate}
$$
\xymatrix{
B^{\oplus n} \ar[r]^{\footnotesize{\begin{pmatrix}0 \\ g\varphi_{(m',n)} \end{pmatrix}}} 
\ar[d]_{gE_n} & B^{\oplus n'}\oplus B^{\oplus m'} 
\ar[d]^{\footnotesize{\begin{pmatrix}gE_{n'} & 0 \\ 0 & E_{m'} \end{pmatrix}}}\\B^{\oplus n} \ar[r]_{\footnotesize{\begin{pmatrix}0 \\ g\varphi_{(m',n)} \end{pmatrix}}} 
\ar[ru]_H & 
B^{\oplus n'}\oplus B^{\oplus m'}.
}
$$
If we establish a theory of homotopy natural transformations, 
then it turns out that $\delta$ induces 
a simplicial homotopy natural transformation $j\mu_1\Rightarrow_{\simp}\eta$ 
(for the definition of simplicial homotopy natural transformations, 
see Definition~\ref{df:simplicial homotopy natural transformation} below). 
Therefore by the theory below, there exists a zig-zag sequence of simplicial natural transformation which connects $\eta$ and $j\mu_1$. 
Thus $\eta$ is homotopic to $0$. 
We complete the proof of assertion $\mathrm{(\beta)}$. 
\qed

\sn
The rest of this subsection, we will establish a theory of homotopy natural 
transformations and justify the argument above.

\begin{conv}
\label{conv:chain complex}
For simplicity we set $\cE=\Ch_b(\cM_B(1))$. 
The functor $C\colon \cE\to\cE$ is given by sending a chain complex 
$x$ in $\cE$ to $Cx:=\Cone\id_x$ 
the canonical mapping cone of 
the identity morphism of $x$. 
Namely the degree $n$ part of $Cx$ is ${(Cx)}_n=x_{n-1}\oplus x_n$ 
and the degree $n$ boundary morphism $d_n^{Cx}\colon {(Cx)}_n\to {(Cx)}_{n-1}$ 
is given by 
$\displaystyle{d^{Cx}_n=
\begin{pmatrix}
-d_{n-1}^x & 0 \\ 
-\id_{x_{n-1}} & d^x_n 
\end{pmatrix}}$. 
For any complex $x$, we define $\iota_x\colon x\to C(x)$ and 
$r_x\colon CC(x) \to C(x)$ to be chain morphisms by setting 
${(\iota_x)}_n=
\begin{pmatrix}
0 \\ \id_{x_n}
\end{pmatrix}$ and 
$\displaystyle{{(r_x)}_n=\begin{pmatrix}
0 & \id_{x_{n-1}} & \id_{x_{n-1}} & 0\\
0 & 0 & 0 &\id_{x_n} 
\end{pmatrix}}$. 

We can show that a pair of chain morphisms $f$, $g\colon x\to y$ in $\cE$ 
are chain homotopic if and only if there exists a morphism 
$H\colon C(x)\to y$ such that $f-g=H\iota_x$. 
We denote this situation by $H\colon F\Rightarrow_C  g$ and we 
say that $H$ is a {\it $C$-homotopy from $f$ to $g$}. 
We can also show that for any complex $x$ in $\cE$, 
$r_x$ is a $C$-homotopy from $\id_{Cx}$ to $0$. 

Let $[f\colon x\to x']$ and $[g\colon y\to y']$ 
be a pair of objects in $\cE^{[1]}$ the morphisms 
category of $\cE$. 
A ($C$-){\it homotopy commutative square} ({\it from $[f\colon x\to x']$ 
to $[g\colon y\to y']$}) 
is a triple $(a,b,H)$ consisting of 
chain morphisms $a\colon x\to y$, $b\colon x'\to y'$ and $H\colon Cx\to y'$ 
in $\cE$ such that $H\iota_x=ga-bf$. 
Namely $H$ is a $C$-homotopy from $ga$ to $bf$. 

Let $[f\colon x\to x']$, $[g\colon y\to y']$ and $[h\colon z\to z']$ 
be a triple of objects in $\cE^{[1]}$ and let 
$(a,b,H)$ and $(a',b',H')$ be homotopy commutative squares from 
$[f\colon x\to x']$ to $[g\colon y\to y']$ and from $[g\colon y\to y']$ 
to $[h\colon z\to z']$ respectively.  
Then we define $(a',b',H')(a,b,H)$ to 
be a homotopy commutative square from 
$[f\colon x\to x']$ to $[h\colon z\to z']$ by setting
\begin{equation}
\label{eq:comp df}
(a',b',H')(a,b,H):=(a'a,b'b,H'\star H)
\end{equation}
where $H'\star H$ is a $C$-homotopy from $ha'a$ to $b'bf$ 
given by the formula
\begin{equation}
\label{eq:star df}
H'\star H:=b'H+H'Ca.
\end{equation}

We define $\cE_h^{[1]}$ to be a category 
whose objects are morphisms in $\cE$ and 
whose morphisms are homotopy commutative squares and 
compositions of morphisms are given by the formula 
$\mathrm{(\ref{eq:comp df})}$ and 
we define $\cE^{[1]}\to\cE_h^{[1]}$ 
to be a functor by sending an object 
$[f\colon x\to x']$ to $[f\colon x\to x']$ 
and a morphism $(a,b)\colon [f\colon x\to x']\to [g\colon y\to y']$ to 
$(a,b,0)\colon [f\colon x\to x']\to [g\colon y\to y']$. 
By this functor, we regard $\cE^{[1]}$ 
as a subcategory of $\cE^{[1]}_h$. 

We define $Y\colon\cE_h^{[1]}\to \cE$ which sending an object 
$[f\colon x\to y]$ to $Y(f):=y\oplus C(x)$ and a homotopy commutative square 
$(a,b,H)\colon [f\colon x\to y]\to [f'\colon x'\to y']$ to 
$\displaystyle{Y(a,b,H):=\begin{pmatrix}b & -H\\ 0 & Ca \end{pmatrix}}$.

We write $s$ and $t$ for the functors $\cE_h^{[1]}\to\cE$ 
which sending an object 
$[f\colon x\to y]$ to $x$ and $y$ respectively. 
We define $j_1\colon s\to Y$ and $j_2\colon t\to Y$ 
to be natural transformations by setting 
$\displaystyle{{j_{1}}_f:=\begin{pmatrix}f\\ -\iota_x \end{pmatrix}}$ and 
$\displaystyle{{j_{2}}_f:=\begin{pmatrix}\id_y\\ 0 \end{pmatrix}}$ 
respectively for any object $[f\colon x\to y]$ in $\cE_h^{[1]}$. 
\end{conv}

\begin{df}
\label{df:homotopy natural transformations}
Let $\cI$ be a category and let $f$, $g\colon \cI\to\cE$ be 
a pair of functors. 
A {\it homotopy natural transformations} ({\it from $f$ to $g$}) 
is consisting of a family of morphisms 
$\{\theta_i\colon f_i\to g_i\}_{i\in\Ob\cI}$ 
indexed by the class of objects of 
$\cI$ and a family of $C$-homotopies 
$\{\theta_a\colon g_a\theta_i\Rightarrow_C \theta_j f_a\}_{a\colon i\to j\in\Mor\cI}$ indexed by the class of morphisms of 
$\cI$ such that for any object $i$ of $\cI$, $\theta_{\id_i}=0$ 
and for any pair of composable morphisms $i\onto{a}j\onto{b}k$ in $\cI$, 
$\theta_{ba}=\theta_b\star\theta_a(=g_b\theta_a+\theta_bCf_a)$. 
We denote this situation by $\theta\colon f\Rightarrow g$. 
For a usual natural transformation $\kappa\colon f\to g$, 
we regard it as a homotopy natural transformation by setting 
$\kappa_a=0$ for any morphism $a\colon i\to j$ in $\cI$. 

Let $h$ and $k$ be another functors from $\cI$ to $\cE$ and 
let $\alpha\colon f\to g$ and $\gamma\colon h\to k$ be natural 
transformations and $\beta\colon g\Rightarrow h$ 
a homotopy natural transformation. 
We define $\beta\alpha \colon f\Rightarrow h$ and 
$\gamma\beta\colon g\Rightarrow k$ to be homotopy natural transformations by 
setting for any object $i$ in $\cI$, 
${(\beta\alpha)}_i=\beta_i\alpha_i$ and ${(\gamma\beta)}_i=\gamma_i\beta_i$ 
and for any morphism $a\colon i\to j$ in $\cI$, 
${(\beta\alpha)}_{a}:=\beta_aC(\alpha_i)$ and 
${(\gamma\beta)}_a=\gamma_j\beta_a$.
\end{df}

\begin{ex}
\label{ex:homotopy natural transformations}
We define $\epsilon\colon s\Rightarrow t$ and $p\colon Y\Rightarrow t$ to be 
homotopy natural transformations between functors 
$\cE_h^{[1]}\to\cE$ by setting for any object 
$[f\colon x\to y]$ in $\cE_h^{[1]}$, 
$\epsilon_f:=f\colon x\to y$ and 
$p_f:=\begin{pmatrix}\id_y & 0\end{pmatrix}\colon Y(f)=y\oplus Cx\to y$ 
and for a homotopy commutative square 
$(a,b,H)\colon [f\colon x\to y]\to [f'\colon x'\to y']$, 
$H\colon f'a\Rightarrow_C bf $ and 
$\displaystyle{p_{(a,b,H)}:=
\begin{pmatrix}0 & -Hr_x \end{pmatrix}\colon b
\begin{pmatrix} \id_y & 0\end{pmatrix}\Rightarrow_C 
\begin{pmatrix} \id_{y'}& 0
\end{pmatrix}
\begin{pmatrix} 
b & -H\\ 0 & Ca
\end{pmatrix}}$.  
Then we have the commutative diagram of homotopy natural transformations.
$$
\xymatrix{
s \ar[r]^{j_1} \ar[rd]_{\epsilon} & Y \ar[d]_p & 
t \ar[l]_{j_2} \ar[ld]^{\id_t}\\
& t.
}
$$
Here we can show that for any object $[f\colon x\to y]$ in $\cE_h^{[1]}$, 
$p_f$ and ${j_2}_f$ are chain homotopy equivalences. 
In particular if $f$ is a chain homotopy equivalence, 
then ${j_1}_f$ is also a chain homotopy equivalence.
\end{ex}

\begin{df}[\bf Mapping cylinder functor on $\Nat_h(\cE^{\cI})$]
Let $\cI$ be a small category. 
We will define $\Nat_h(\cE^{\cI})$ {\it the category of homotopy natural 
transformations} ({\it between the functors from $\cI$ to $\cE$}) 
as follows. 
An object in $\Nat_h(\cE^{\cI})$ is a triple $(f,g,\theta)$ 
consisting of functors $f$, $g\colon \cI\to\cE$ and a homotopy 
natural transformation $\theta\colon f\Rightarrow g$. 
A morphism $(a,b)\colon(f,g,\theta)\to(f',g',\theta')$ 
is a pair of natural transformations 
$a\colon f\to f'$ and $b\colon g\to g'$ such that $\theta'a=b\theta$. 
Compositions of morphisms is given by component-wise compositions of natural 
transformations. 

We will define functors $S$, $T$, $Y\colon\Nat_h(\cE^{\cI})\to\cE^{\cI}$ 
and natural transformations $J_1\colon S\to Y$ and 
$J_2\colon T\to Y$ as follows. 
For any object $(f,g,\theta)$ and any morphism 
$(\alpha,\beta)\colon (f,g,\theta)\to (f',g',\theta')$ in 
$\Nat_h(\cE^{\cI})$ 
and any object $i$ and 
any morphism $a\colon i\to j$ in $\cI$, we set
$$S(f,g,\theta):=f,\ S(\alpha,\beta):=\alpha,$$
$$T(f,g,\theta):=g,\ T(\alpha,\beta):=\beta,$$
$${Y(f,g,\theta)}_i(={Y(\theta)}_i):=Y(\theta_i)(=g_i\oplus C(f_i)),$$
$${Y(f,g,\theta)}_a(={Y(\theta)}_a):=Y(f_a,g_a,\theta_a)\left(  
\begin{pmatrix}g_a & -\theta_a \\ 0 & Cf_a \end{pmatrix}
\right ),\  
{Y(\alpha,\beta)}_i:=
\begin{pmatrix}\beta_i & 0 \\ 0 & C(\alpha_i) \end{pmatrix},$$
$${J_1}_{{(f,g,\theta)}_i}:={j_1}_{\theta_i},\ 
{J_2}_{{(f,g,\theta)}_i}:={j_2}_{\theta_i}$$
In particular for an object $(f,g,\theta)$ 
in $\Nat_h(\cE^{\cI})$ if for any object 
$i$ of $\cI$, $\theta_i$ is a chain homotopy equivalence, 
then there exists a zig-zag sequence of morphisms which connects $f$ to 
$g$, $f\onto{J_1}Y(\theta)\overset{J_2}{\leftarrow}g$ such that for any object $i$, ${J_1}_i$ and ${J_2}_i$ are 
chain homotopy equivalences. 
\end{df}

\begin{df}[\bf Simplicial homotopy natural transformation]
\label{df:simplicial homotopy natural transformation}
Let $\cI$ be a simplicial small category, namely a simplicial object in 
the category of small categories. 
We write $\Ch_b(S_{\cdot}\cM_B(1))$ for 
the simplicial additive category defined by 
sending $[n]$ to $\Ch_b(S_n\cM_B(1))$ 
the category of bounded chain complexes on $S_n\cM_B(1)$. 
Notice that there is a canonical identification $S_{\cdot}\Ch_b(\cM_B(1))=\Ch_b(S_{\cdot}\cM_B(1))$. 
Let $f$, $g\colon\cI\to\Ch_b(S_{\cdot}\cM_B(1))$ be simplicial functors. 
Recall that a {\it simplicial natural transformation} ({\it from $f$ to $g$}) 
is a family of natural transformations $\{\rho_n\colon f_n\to g_n\}_{n\geq 0}$ 
indexed by non-negative integers such that for any morphism 
$\varphi\colon [n]\to [m]$, 
we have the equality $\rho_nf_{\varphi}=g_{\varphi}\rho_m$. 

A {\it simplicial homotopy natural transformation} ({\it from $f$ to $g$}) 
is a family of homotopy natural transformations $\{\theta_n\colon f_n\Rightarrow g_n\}_{h\geq 0}$ 
indexed by non-negative integers 
such that for any morphism 
$\varphi\colon [n]\to [m]$, 
we have the equality $\theta_nf_{\varphi}=g_{\varphi}\theta_m$. 
We denote this situation by $\theta\colon f\Rightarrow_{\simp} g$. 

For a simplicial homotopy natural transformation 
$\theta\colon f\Rightarrow_{\simp} g$, 
we will define $\cY(\theta)\colon \cI\to\Ch_b(S_{\cdot}\cM_B(1))$ and 
$\fJ_1\colon f\to \cY(\theta)$ and $\fJ_2\colon g\to \cY(\theta)$ 
to be a simplicial functor and simplicial natural transformations respectively 
as follows. 
For any $[n]$ and any morphism $\varphi\colon[m]\to[n]$, 
we set ${\cY(\theta)}_n=Y(\theta_n)$, ${\fJ_1}_n:={J_1}_{\theta_n}$, 
${\fJ_2}_n:={J_2}_{\theta_n}$ and ${\fY(\theta)}_{\varphi}:=Y(f_{\varphi},g_{\varphi})$. 
In particular if for any non-negative integer $n$, any object $j$ of $\cI_n$, 
${\theta_n}_j$ is a chain homotopy equivalence, then there exists a zig-zag sequence of morphisms which connects $f$ to $g$, $f\onto{\fJ_1}\cY(\theta)\overset{\fJ_2}{\leftarrow}g$ such that for any non-negative integer $n$ and any object $j$, ${{\fJ_1}_n}_j$ and ${{\fJ_2}_n}_j$ are chain homotopy equivalences.
\end{df}

\paragraph{Acknowledgement}
I wish to express my deep gratitude to Seidai Yasuda, 
Marco Schlichting, Osamu Iyama and the referees for carefully reading 
a version of this paper and giving me valuable comments to make the paper 
more readable. 
I also greatful for Kazuya Kato for teaching me deeply impressive suggestions. 
I also thank to Norihiko Minami for inviting me to a workshop in a memory of 
Tetsusuke Ohkawa and giving me the opportunity to talk about the 
contents of this article (the slide movie of my talk at the workshop is \cite{Moc15}).

\mn
SATOSHI MOCHIZUKI\\
\emph{DEPARTMENT OF MATHEMATICS,
CHUO UNIVERSITY,
BUNKYO-KU, TOKYO, JAPAN.}\\
e-mail: {\tt{mochi@gug.math.chuo-u.ac.jp}}

\end{document}